\input amstex
\documentstyle{amsppt}
\magnification=\magstep1
\TagsOnRight

\hoffset1 true pc
\voffset2 true pc
\hsize36 true pc

\vsize52 true pc

\tolerance=2000
\define\m1{^{-1}}
\define\ov1{\overline}
\def\gp#1{\langle#1\rangle}
\def\cry#1{\operatorname{\frak{Crys}}(#1)}
\def\ker#1{\operatorname{{ker}}#1}

\catcode`\@=11
\def\logo@{}
\catcode`\@=\active

\topmatter
\title
Torsion free groups with indecomposable holonomy group I
\endtitle
\author
V.A.~Bovdi, P.M.~Gudivok, V.P.~Rudko
\endauthor
\dedicatory
Dedicated to Professor {\it L.G.~Kov\'acs}  on his 65th birthday
\enddedicatory
\leftheadtext\nofrills{ V.A.~Bovdi, P.M.~Gudivok, V.P.~Rudko }
\rightheadtext \nofrills { Torsion free groups with indecomposable
holonomy group } \abstract We study  the  torsion free generalized
crystallographic groups with the indecomposable holonomy group
which is isomorphic to either a cyclic group of order ${p^s}$ or a
direct product of two cyclic groups of order ${p}$.
\endabstract
\thanks
The research was supported by OTKA  No.T 025029
\endthanks
\address
\hskip-\parindent
{\sl  V.A.~Bovdi}
\newline
Institute of Mathematics and Informatics
\newline University of Debrecen
\newline
H-4010 Debrecen, P.O.  Box 12
\newline Hungary
\newline
vbovdi\@math.klte.hu
\bigskip
\hskip-\parindent
{\sl P.M.~Gudivok, V.P.~Rudko}
\newline
Department of Algebra
\newline Uzhgorod  University
\newline
88 000, Uzhgorod
\newline Ukraina
\newline
math1\@univ.uzhgorod.ua
\endaddress\endtopmatter
\document
\subhead
 Introduction
\endsubhead
A classical {\it crystallographic group} is a discrete cocompact
subgroup of $I(\Bbb R^m)$, the isometry group of $\Bbb R^m$.
Torsion-free crystallographic groups are called {\it Bieberbach
groups}. The present state of the theory of crystallographic
groups and a historical overview as well as its connections to
other branches of mathematics are described in \cite{16,17}.

In this paper we consider generalized torsion-free
crystallographic groups  with indecomposable holonomy groups
isomorphic to either $C_{p^s}$ or $C_p\times C_p$.

It was shown in  \cite{6, 7, 13} that the description of the
$n$-dimensional crystallographic groups for an arbitrary $n$ is of
wild  type, in the sense that it relates to the classical
unsolvable problem of describing the canonical forms of pairs of
linear operators acting on finite dimensional vector spaces.

Using  Diederichsen's  classification of integral representations
of the cyclic group of prime order (see \cite{5}), L.~Charlap
\cite{4} gave a full classification of Bieberbach groups with the
cyclic  holonomy group $G$ of prime order. G.~Hiss and
A.~Szczep\'anski   \cite{12}  proved that Bieberbach groups with
nontrivial irreducible holonomy group $G$ do not exist. G.~Kopcha
and V.~Rudko  \cite{13}  studied    torsion-free crystallographic
groups with indecomposable cyclic holonomy group of order $p^n$
the classification of which for $n\geq 5$ also has wild type.

J.~Rossetti and P.~Tirao \cite{18, 19, 20} described the
torsion-free crystallographic groups whose holonomy group are
direct sums of indecomposable subgroups of $GL(n, \Bbb Z)$ ($n\leq
5$) and isomorphic to $C_2\times C_2$.

Similarly, interesting results on this topic were obtained in the
research of N.~Gupta and S.~Sidki  \cite{8, 9}.

We need the following definitions and notation for announcing our
results.

Let $K$ be a principal domain, $F$ be a field containing $K$ and
let $G$ be a finite group.  Let $M$ be a $KG$-module of a faithful
matrix $K$-representation $\Gamma$ of $G$ and let $FM$ be a linear
space over $F$ in which the $K$-module $M$ is a full lattice.
Let $\widehat{M}=FM^+/M^+$ be the quotient group of the additive group
$FM^+$ of the linear space $FM$ by the additive group $M^{+}$ of
the module $M$. Then $FM$ is an  $FG$-module and $\widehat{M}$ is a
$KG$-module with operations:
$$
g\cdot (\alpha m)=\alpha g(m); \,\,\,\,\, g\cdot (x+M)=g(x)+M,
$$
where $g\in G$, $\alpha\in F$, $m\in M$,
$x\in FM$.

Let $T:G\to\widehat{M}$ be a $1$-cocycle  of $G$ with values in
$\widehat{M}$; that is, $T(g)$ is regarded as the set $x_1+M$,
where $x_1$ is an element from the coset of $T(g)$ in
$\widehat{M}$. We define the group
$$
\cry{G; M; T}=\{\,\,\,\, (g,x)\,\,\,\, \mid\,\,\,
g\in G, \,\,\, x\in T(g)\,\,\,\, \}
$$
with the operation$$(\,g,\,x\,)\cdot (\,g',\,x'\,)=(\,gg',\,
g'x+x'\,),$$ where $g,g'\in G$, $x\in T(g)$, $x'\in T(g')$.

The purpose of  this paper is to study the group $\cry{G; M; T}$,
and in particular to determine  when  $\cry{G; M; T}$ is a torsion
free group. We note  that if $K=\Bbb Z$ and $F=\Bbb R$, then
$\cry{G; M; T}$ is isomorphic to an  $n$-dimensional classical
crystallographic group, where $n=rank_{\Bbb Z}M$.

We use the terminology of
the theory of group  representations.   The group $\cry{G; M; T}$
is called {\it irreducible} ( {\it indecomposable} ), if $M$ is an
irreducible ( indecomposable ) $KG$-module and the cocycle $T$ is
not cohomologous  to  zero.

A cocycle $T:G\to \widehat{M}$ is called {\it coboundary}, if
for every $ g\in G$ there exists an $x\in FM$ such  that $T(g)=(g-1)x+M$.
The cocycles $T_1:G\to \widehat{M}$ and $T_2: G\to \widehat{M}$ are
called {\it cohomologous} if $T_1-T_2$ is a coboundary.

Let $C^{1}(G, \widehat{M})$, $B^{1}(G, \widehat{M})$,
$H^{1}(G,\widehat{M}) =C^{1}(G, \widehat{M})/B^{1}(G,
\widehat{M})$ be the group of the cocycles, the coboundaries and
the cohomologies  of $G$ with values in the group $\widehat{M}$,
respectively. The group $\cry{G; M; T}$ is an extension of $M^{+}$
by  the group $G$.   This extension splits if and only if
$T\in B^{1}(G,\widehat{M})$. Therefore, the group $\cry{G; M; T}$
splits for all $T$ if and only if $H^{1}(G,\widehat{M})$ is  trivial.

\subhead
 Main Results
\endsubhead
Using  results from  \cite{2, 3,  10, 11, 14, 15}, we prove the
following three theorems. We note that Lemma 12 can considered as
a separate result.

\proclaim {Theorem 1 } Let $K$ be either the ring of rational
integers $\Bbb Z$, or $p$-adic integers $\Bbb  Z_{p}$, or the
localization  $\Bbb Z_{(p)}$  of $\Bbb Z$ at $p$
and let $G\cong C_{p^s}$ be a cyclic group of order $p^s$.   If
$s\geq 3$, then the set of $K$-dimensions of the indecomposable
$KC_{p^s}$-modules $M$, for which there exist  torsion-free groups
$\cry{C_{p^s}; M; T}$, is unbounded.
\endproclaim

\proclaim {Theorem 2 } Let $K$ be either the localization $\Bbb
Z_{(p)}$  of $\Bbb Z$ at $p$,   or the ring of $p$-adic integers
$\Bbb Z_{p}$  and let $G=\gp{a}\cong C_{p^2}$ be a
cyclic group of order $p^2$.  Up to isomorphism, all torsion-free
indecomposable groups $\cry{C_{p^2}; M;  T}$ can be described
with the help of the following indecomposable $KC_{p^2}$-modules
$M$ and cocycles $T$ of the group $C_{p^2}$ with values in the
groups $\widehat{M}=FM^+/M^{+}$:
\itemitem{1) } $M=X_i=\gp{\,\, (a-1)\Phi(a^p),\,\,\,
\Phi(a)+(a-1)^{i+1}\,\,\, }$,
$T=T_i$, where
$\Phi(x)=x^{p-1}+\cdots+x+1$,
$T_i(a)=p^{-2}\Phi(a)\Phi(a^p)+X_i$,
and $i=0,1,\ldots,p-2$;
\itemitem{2)} $p>2$ and
$M=U_j=\gp{\,\,\, \big( (a-1)^{j+1}+\Phi(a), \,\,\,
(a-1)^{j}\big);\,\,\, \Phi(a^p)(a-1,1) \,\,\,}$, a
$KC_{p^2}$-submodule in $ \big(KC_{p^2}\big)^{(2)}=\{ \,\,
(x_1,x_2) \,\, \mid \,\, x_1,x_2\in KC_{p^2} \,\, \} $, and
$T=f_j$, where $f_j(a)=p^{-2}\Phi(a)\Phi(a^p)(1,0)+U_j$ and
$j=1,\ldots,p-2$.

The number of these groups $\cry{C_{p^2};  M;  T}$ is equal to
$2p-3$.
\endproclaim

\proclaim
{Corollary 1}
There  exist at least $2p-3$ Bieberbach ( in the classical sense)
groups having a cyclic indecomposable holonomy group of order ${p^2}$.
\endproclaim

\proclaim {Theorem 3} Let $G=\gp{a}\times \gp{b}\cong C_2\times
C_2$ and let $K$ be either the  ring of rational integers $\Bbb
Z$,  or the ring of $2$-adic integers $\Bbb Z_{2}$,  or the
localization  $\Bbb Z_{(2)}$  of $\Bbb Z$ at the prime  $2$.
Furthermore, let  $F$ be a field containing $K$,  and $M$ be a $KG$-module
of the indecomposable $K$-representation $\Gamma$ of $G$, and let
$\cry{G; M; f}$ be the group defined by the cocycle $f:G\to
\widehat{M}=FM^{+}/M^{+}$. The following table exhibits
indecomposable $K$-representations $\Gamma$ of $G$ and
cocycles $f$, which define, up to isomorphism, all torsion-free
indecomposable groups  $\cry{G; M; f}$:

{ \eightpoint{
\centerline{\vbox{\halign{\strut\offinterlineskip\vrule\quad\hfill
$#$\hfill&\quad
\vrule\quad\hfill$#$\hfill&\quad\vrule\quad\hfill$#$\hfill&\quad
\vrule\quad\hfill$#$\hfill&\quad\vrule\quad\hfill$#$\hfill\quad\vrule\cr
\noalign{\hrule} &&&&\cr \noalign{\vskip -5pt}
N:&m&\Gamma&f(a)=(x_1,\ldots, x_m)+M,\,\,\, f(b)=(y_1,\ldots,
y_m)+M&t_m\cr \noalign{\vskip -3pt} &&&&\cr \noalign{\hrule}
&&&&\cr \noalign{\vskip -3pt} 1& 4n+1& \Delta_n&
x_{n+1}=\frac{1}{2},\,\,\,\,\,  x_i=0  \,\,\,\,\, (i\not=n+1), &
2^{n-1}\cr \noalign{\vskip -1pt} &(n\geq 1)&
&y_1=\frac{1}{2},\,\,\,\,\,   2y_{2}=\cdots=2y_{n+1}=0,  &\cr
\noalign{\vskip -1pt} &&&y_2+\cdots+y_{n+1}=\frac{1}{2},&\cr
\noalign{\vskip -1pt} &&&y_{n+2}=\cdots=y_{4n+1}=0&\cr
\noalign{\vskip -3pt} &&&&\cr \noalign{\hrule} &&&&\cr
\noalign{\vskip -3pt} 2
& 4n+4                                           & W_n^*
& x_{2n+3}=\frac{1}{2},\,\,\,\,  x_i=0 \,\,\,\, (i\not=2n+3),
& 2^n \cr\noalign{\vskip -1pt}
                                                &(n\geq 0)                                       &                                                 & y_1=0, \,\,\,\,\,\, y_2=\frac{1}{2},
\,\,\,\,\,\, y_3=\cdots=y_{3n+3}=0,                          & \cr
\noalign{\vskip -1pt} & & & 2y_{3n+4}=\cdots=2y_{4n+3}=0,
\,\,\,\,\,\, y_{4n+4}=\frac{1}{2} & \cr \noalign{\vskip -3pt}
&&&&\cr\noalign{\hrule} &&&&\cr \noalign{\vskip -3pt} 3
& {5}                                      & \Delta_1^* &
f(a)=(0,\frac{1}{2},0,0,0), \,\,\,\,\,\,
f(b)=(\frac{1}{2},0,\frac{1}{2},\frac{1}{4},0)&1 \cr
\noalign{\vskip -3pt} &&&&\cr \noalign{\hrule} &&&&\cr
\noalign{\vskip -3pt} 4 &8 & W_1 &
f(a)=(0,0,0,0,\frac{1}{2},0,0,0), \,\,\,\,\,\,
f(b)=(0,\frac{1}{2},0,\frac{1}{4},0,\frac{1}{2},0,\frac{1}{2})  &
1 \cr \noalign{\vskip -3pt} &&&&\cr \noalign{\hrule}}}} } }
\noindent where $m$ is the degree of the  representation $\Gamma$,
$f$ is  the cocycle, and $t_m$ is the number of the torsion free
groups.
\endproclaim

\subhead
  Preliminary results
\endsubhead
Let $K=\Bbb Z$, ${}\Bbb  Z_{(p)}$ or $\Bbb  Z_{p}$ as above.  We
point out that in these cases  the group $H^1(G,\widehat{M})$ is
finite. Denote by  $C_{p^n}=\gp{a\mid a^{p^n}=1 }$ the cyclic
group of order $p^n$. The following three Lemmas and Corollary 2
are well-known and they can be found for example in \cite{1}.

\proclaim {Lemma 1} Let $K$ be either   the ring of rational
integers $\Bbb Z$,  or the ring of $p$-adic integers $\Bbb Z_{p}$,
or the localization  $\Bbb Z_{(p)}$  of $\Bbb Z$ at $p$,
respectively. Let $G_i$ ($i=1,2$) be a finite group and
$\Gamma_i$, $M_i$, $T_i$ ($i=1,2$) be
the representation, the module and the cocycle associated with
$G_i$ as in the introduction. The groups $\cry{G_1; M_1;
T_1}$ and $\cry{G_2; M_2; T_2}$ are isomorphic if and only if
there exist a group isomorphism $\varepsilon: G_1\to G_2$ and a
$K$-module isomorphism $\tau:M_1\to M_2$ which satisfy the
following conditions:
\itemitem{1)}  $\varepsilon(g)\tau=\tau g$  in $M_1$, for all
$g\in G_1$;
\itemitem{2)} the cocycles $T_2$ and $T_1^\varepsilon$ are
cohomologous (here,  $T_1^\varepsilon(g)=\tau' T_1(\varepsilon \m1g)$,
$g\in G_2$, where $\tau':\widehat {M_1}\to \widehat {M_2}$ is the
homomorphism induced  by the  homomorphism $\tau$).
\endproclaim

\rightline{\text{\qed}}

\proclaim {Lemma 2 } Assume that the character of the
$K$-representation $\Gamma$ of the group $C_n$ does not  contain
the trivial character as a summand.  Then
$H^1(C_n,\widehat{M})$ is  trivial.
\endproclaim

\demo{Proof}
Since $1$ is not an eigenvalue of the
 operator $a$, which acts on
$FM$,  the operator $a-1$ is a unit.  This means that
$T(a)=(a-1)x+M$ for some $x\in FM$, i.e.
$B^{1}(C_n,\widehat{M})=C^1(C_n,\widehat{M})$.
\enddemo
\rightline{\text{\qed}}

\proclaim
{Lemma 3}
Let $G\cong C_{p^s}$ and $M$ be a projective $KG$-module.
Then $H^1(C_{p^s},\widehat{M})$ is trivial.
\endproclaim

\demo{Proof} It is well known that  multiple direct sums
$M\oplus\cdots\oplus M$ of the module $M$ are  free
$KC_{p^s}$-modules. Therefore, it is sufficient to prove the
lemma for $M=KC_{p^s}$. Let $T(a)=x+M$ for some
$$
x=\lambda(1+a+\cdots+a^{p^s-1})+u_1(a-1)\in FM,
$$
where $\lambda\in F$ and $u_1\in FC_{p^s}$.  From
the condition $(1+a+\cdots+a^{p^s-1})T(a)\subset M$ it
follows that $\lambda p^s\in K$.   Then $x-\lambda
p^s=u_2(a-1)$, where $u_2 \in FC_{p^s}$.
Therefore, $T$ is  coboundary.
\enddemo\rightline{\text{\qed}}

\proclaim {Corollary 2}
Assume that a  $K$-representation
$\Gamma$ of the group $C_p$ does not contain  the trivial
$K$-representation as a summand.   Then $H^1(C_p,\widehat{M})$ is
trivial.
\endproclaim
\demo{Proof}
A $K$-representation $\Gamma$ of the group $C_p$ is a
direct sum $\Gamma=\Gamma_1\oplus\Gamma_2$, where
$\Gamma_1$ is a sum of the  irreducible $K$-representation of
degree $p-1$,  and $\Gamma_2$ is a  $K$-representation
corresponding   to a projective $KC_p$-module. The proof follows by
applying Lemma 1 to $\Gamma_1$ and  Lemma 3 to $\Gamma_2$.
\enddemo\rightline{\text{\qed}}

For the proof of  Theorem 1  we will consider a $K$-representation
of the group $\gp{a}\cong C_{p^s}$.   Let $\xi_t$ be a primitive
$p^t$th root of unity and $\xi_{t-1}=\xi_t^p$ ($t\geq 1$).  Put
$$
B_0=\{1\},\,\,\,\,
B_1=\{1,\xi_1,\ldots,\xi_1^{p-2}\},
\,\,\,\,\,
B_j=\cup_{i=0}^{p-1}\xi_j^iB_{j-1}
\,\,\,\,\, (j\geq 2).
$$
 Thus $B_t$ is a $K$-basis of the
ring $K_t=K[\xi_t]$, which is a $KC_{p^s}$-module: $
a(\alpha)=\xi_t\alpha, $ where $\alpha\in K_t$, $t\leq s$. The set
$B_t$ is  also be an $F$-basis of the space $FK_t$ ($t\geq 0$).

Let $\delta_t$ be a $K$-representation of $C_{p^s}$
corresponding to the $K$-basis of the module $K_t$ of this
representation.   We note that $\delta_t$ is an irreducible
$K$-representation of $C_{p^s}$ and
$\delta_t^p(a)=E_p\otimes \delta_{t-1}(a)$, where
$E_p$ is the  identity matrix of degree $p$.
 Let
$$
\Delta_1=\delta_0^{(n)}+\delta_1^{(n)};
\,\,\,\,\,\,\,\,\,
\Delta_2=\delta_2^{(n)}+\delta_s^{(n)}
$$
be the  sums  of $2n$ irreducible $K$-representations of
$C_{p^s}$, where
$\delta_i^{(n)}=\underbrace{\delta_i+\cdots+\delta_i}_n$. Let us
consider the following $K$-representation $\Delta$ of the group
$C_{p^s}$:
$$
\Delta(a)=
\left(\smallmatrix
\Delta_1(a) &  U(a)      \\
   0        & \Delta_2(a)\\
\endsmallmatrix\right),$$
where
$$
U(a)=
\left(\smallmatrix
E_n\otimes\gp{1}_0  & J_n\otimes\gp{1}_0        \\
E_n\otimes\gp{1}_1  & J_n\otimes\gp{1}_1        \\
\endsmallmatrix\right);$$

$J_n$ is the Jordan block of  order $n$ with ones on the main
diagonal and the symbol $\gp{\omega}_t$
denotes the matrix with zero columns except the last one, which
consists of the coordinates of the element $\omega\in K_t$ written in the
basis $B_t$ ( $t=0,1$). \rightline{\text{\qed}}

\proclaim
{Lemma 4} ( see \cite{2, 3} )
The $K$-representation $\Delta$ of $C_{p^s}$ is indecomposable.
\endproclaim
\rightline{\text{\qed}}

\proclaim {Lemma 5} Let $x\in FK_t$ ($t>0$), such that $(a-1)x\in
K_t$.   Then $px\in K_t$ and all coordinates of the vector $px$
are multiples  of the last coordinate.
\endproclaim
\demo{Proof} The $K$-basis $B_t$ in $K_t$ will be an $F$-basis in
 $FK_t$.   Let us consider the coordinates of the  column vectors
in  $FK_t$ and  the matrix $\delta_t(a)$ of the operator $a$ in
this basis. The lemma is easily checked successively for $t=1$,
$t=2$ and so on.
\enddemo\rightline{\text{\qed}}


 Let $B$ be a $K$-basis  of the $K$-module
$M_\Delta$ affording the matrix $K$-representation $\Delta$ of the
group $G=\gp{a}\cong C_{p^s}$. Denote   the first basis element by
$v$. It is easy to see that $B$ is an  $F$-basis in $FM_\Delta$.
We  define the following function:
$$
T_\Delta:C_{p^s}\to \widehat{M_\Delta}=FM_\Delta^{+}/M_\Delta^{+},
$$
where $T_\Delta(a^j)=jp^{-s}v+M_\Delta$ and $j=0,1,\ldots,p^s-1$.

\proclaim {Lemma 6} The function $T_\Delta$ is a $1$-cocycle of
the group $C_{p^s}$ with values in the group $\widehat{M_\Delta}$.
The cocycle $T_\Delta$ is not cohomologous to the zero cocycle at
the element $b=a^{p^{s-1}}$ of  order $p$.
\endproclaim
\demo{Proof} The first statement about the  function of $T_\Delta$
follows from the definition of this function.  For proving the
second part of the lemma, let us consider the $p$th power
$\Delta^p(a)$ of the representation $\Delta$.   We note that
$\Delta^p(a)=\left(\smallmatrix
\Delta_1^p(a) &  U'(a)      \\
   0        & \Delta_2^p(a)\\
\endsmallmatrix\right)$ and  $\Delta_1^p(a)=E$.
Clearly,  the first row in $U'(a)$ has the form:
$$
(\,\,\, \gp{1}_0,\ldots, \gp{1}_0,
\gp{1}_0,\ldots,  \gp{1}_0\,\,\,)
$$
and the row  of matrices corresponding to the first of the
representations $\delta_1^p$ will take  the form of the
following matrix:
$$
(\,\,\, \gp{1}_1, \gp{\xi_1}_1,
\ldots, \gp{\xi_1^{p-1}}_1,
\gp{1}_1,\ldots,  \gp{\xi_1^{p-1}}_1\,\,\, ).
$$

If we sum the rows of this matrix and substract the result from
the 1st row in $U'(a)$, then we will get a row in which all the
elements will be multiples to $p$.  So this transformation of rows
in $U'(a)$ corresponds to the replacement  of some base elements
$u\in B$ ($u\not=v$) to $u'=u\pm v$.  Let us carry out this
replacement;  and let $\Delta'$ be a $K$-representation of the
group $C_{p^s}$ in the new $K$-basis of the module $M_\Delta$.  It
is easy to see that the replacement  of the basis did not
influence the values of the function $T_\Delta$.

Let $H=\gp{b \mid b=a^{p^{s-1}}}$ and $\Delta_H'$ be a
restriction of the representation  $\Delta'$ to the
subgroup $H$.   Then
$$
\Delta_H'(b)=
\left(\smallmatrix
\delta_0^{(m_1)}(a) &       U''(b)            \\
         0          &  \delta_1^{(m_2)}(a)    \\
\endsmallmatrix\right),$$
where, as shown  above, all elements of the first row in $U''(b)$
are multiples of $p$.
Let
$M_\Delta=M_1\oplus M_2$ be a decomposition  of the $K$-module
$M_\Delta$ into the direct sum of $M_1$ and $M_2$, corresponding
to the representations $\delta_0^{(m_1)}$ and $\delta_1^{(m_2)}$.

Suppose that the cocycle $T_\Delta$ is cohomologous to the
trivial cocycle at $H$.  Then there exists a vector $x\in
FM_\Delta$,  such that
$$
T_\Delta(b)=(b-1)x+ M_\Delta.
$$
Let $x=x_1+x_2$ ($x_i\in FM_i, i=1,2$).   Since the projection of
$T_\Delta(b)$ on  $FM_2$  is equal to zero (modulo $M_\Delta$), the
projection of $(b-1)x=(b-1)x_2$ on  $FM_2$ also equals to zero.
From Lemma 5 it follows that $px_2\in M_2$.  Let $\lambda$ be the
coefficient of the basis vector $v$ in $(b-1)x$.

  It is easy
to see that $\lambda$  is a sum of products  of the
elements of the first row in $U''(b)$ (these elements are
multiple  to $p$) on the column, which consist of  coordinates
of the vector $x_2$.
From the condition $px_2\in M_2$ it follows that
$\lambda\in K$.  Since  $T_\Delta(b)=p\m1v+M_\Delta$, we have
 $\lambda=p\m1$. But $p\m1$  does not belong  to $K$,
and, therefore,  $\lambda\not\in K$. This
contradiction proves that $T_{\Delta}$ is not cohomologous to  zero
at $H$.
The lemma is proved.
\enddemo\rightline{\text{\qed}}

\subhead
  Proof of  Theorem 1
\endsubhead
Let us consider the group
$\cry{C_{p^s}; M_\Delta; T_\Delta}$.   If  there  exists an
element of prime order in this group, then this order can
only be $p$ and, moreover, the cocycle $T_\Delta$  must be
cohomologous to  the zero cocycle  in the unique  element
$b=a^{p^{s-1}}$ in the group $C_{p^s}$ of prime order $p$.
According to  Lemma 6, this is impossible.  Therefore, the
group $\cry{C_{p^s};  M_\Delta;  T_\Delta}$ has  no torsion
elements.   Moreover,  this group is  indecomposable (see
Lemma 4).
\rightline{\text{\qed}}

\bigskip

For the cyclic group $\gp{a}\cong C_{p^2}$ we wont to find all
those groups $\cry{C_{p^2}; M; T}$ which are torsion free. Put
$$
\Phi(x)=x^{p-1}+x^{p-2}+\cdots+x+1.
$$
There exists  a unit $\theta$ in $KC_{p^2}$ such that
$$
(a-1)^p\Phi(a^p)=p(a-1)\theta \Phi(a^p).
$$
Let $X_i$ be a $KC_{p^2}$-submodule in $KC_{p^2}$,
generated by the following elements:
$$
u=\Phi(a)\Phi(a^p),
\,\,\,\,\,\,\,\,
\omega=(a-1)\Phi(a^p),
\,\,\,\,\,\,\,\,
v=\Phi(a)+(a-1)^{i+1},
$$
where $0\leq i\leq p-2$. It is easy to see that
$$
(a-1)u=0,
\,\,\,\,\,\,\,\,\,
\Phi(a)\omega=0,
\,\,\,\,\,\,\,\,\,
\Phi(a^p)v=u+(a-1)^{i}\omega.
$$
From these equations it follows that the $K$-representation
$\Gamma_i$ of the group $C_{p^2}$ in the $K$-basis
$$
{ \eightpoint{ \centerline{
\vbox{\halign{\strut\offinterlineskip\quad\hfill $#$\hfill&\quad
\quad\hfill$#$\hfill&\quad \quad\hfill$#$\hfill&\quad
\quad\hfill$#$\hfill&\quad \quad\hfill$#$\hfill&\quad
\quad\hfill$#$\hfill\quad\cr
u;&\omega,& a\omega,  & a^2\omega,  & \ldots,  & a^{p-2}\omega;\cr
  &a^{l}v,    & a^{l+p}v, & a^{l+2p}v,  & \ldots,  & a^{l+p(p-2)}v,\cr
}}}}}
$$
($l=0,1,\ldots, p-1$) corresponding to the module $X_i$, has the
following form:
$$
\Gamma_i(a)=
\left(\smallmatrix
1 & 0                   &  \gp{1}_0             \\
  & \delta_1(a)         &  \gp{\alpha_i}_1      \\
  &                     &  \delta_2(a)          \\
\endsmallmatrix\right),
$$
where $\alpha_i=(\xi_1-1)^i$ and $i=0,1,\ldots,p-2$.

\proclaim {Lemma 7} Let  $H=\gp{b\mid b=a^p}$.  The $KH$-module
$X_i\mid_H$ is a direct sum of two $KH$-submodules, one of which
coincides with $Ku$.
\endproclaim
\demo{Proof}
Let us examine  the $K$-submodule $X_i'$ in $X_i$ generated
by the following system of $p^2-1$ elements from $X_i$:
$$
V=\{v, bv,\ldots,b^{p-2}v\},
\,\,\,
(a-1)V, \ldots, (a-1)^{p-2}V,
\,\,\,
v'=(a-1)^{p-1}v+\theta u,
$$$$
bv', \ldots,b^{p-2}v',
\,\,\,
\theta\omega_1, \ldots, \theta\omega_i,
\,\,\,
u+\omega_{i+1},\omega_{i+2}, \ldots, \omega_{p-1},
$$
where $\omega_{j}=(a-1)^j\Phi(b)=(a-1)^{j-1}\omega$ and
$j=1,\ldots,p-1$.

Clear,  $X_i=Ku\oplus X_i'$ is a direct sum of $K$-modules $Ku$
and $X_i'$. To prove the lemma it is sufficient to show that
$X_i'$ is a $KH$-module.   We have:
$$
\matrix
\format  \l\quad        & \l                  & \l    & \c     \\
\Phi(b)v                &=u+\omega_{i+1}      & \in   & X_i',  \\
\Phi(b)(a-1)v           &=\omega_{i+2}        & \in   & X_i',  \\
\cdots\cdots\cdots\cdots      & \cdots\cdots        & & \cdots \\
\cdots\cdots\cdots\cdots      &\cdots \cdots        & & \cdots \\
\Phi(b)(a-1)^{p-i-2}v   &=\omega_{p-1}        & \in   & X_i',  \\
\Phi(b)(a-1)^{p-r-2+j}v &=p\theta\omega_{j}   & \in   & X_r',  \\
\endmatrix
$$
where $0<i$,\quad  $0\leq r\leq p-2$, \quad  $j=1, \ldots, r$ and
$$
\split
\Phi(b)v'&=(a-1)^{p-1}\Phi(b)v +p\theta u=
(a-1)^{p+i}\Phi(b)+p\theta u\\
&=p\theta (a-1)^{1+i}\Phi(b)+p\theta u=
p\theta(\omega_{i+1}+ u)\in X_i'.
\endsplit
$$

These equations  show that $X_i'$ is the $KH$-submodule
in $X_i$.
\enddemo\rightline{\text{\qed}}

Let us introduce the cocycle
$$
T_i:C_{p^2}\to \widehat{X_i}=FX_i^{+}/X^{+}_i,\tag1
$$
where $T_i(a)=p^{-2}u+X_i$ and $i=0,1,\ldots, p-2$.

\proclaim
{Lemma 8}
The group  $\cry{C_{p^2}; X_i; T_i}$ is  torsion free
($i=0,1,\ldots,p-2$).
\endproclaim
\demo{Proof} Since $T_i(a^p)=p\m1u+X_i\not=X_i$,  from Lemma 7 it
follows that
$$
\big((a^p-1)FX_i+X_i\big)\cap \big(Fu+X_i\big)=X_i.
$$
These conditions  show  that the cocycle $T_i$ is  not
cohomologous to  the zero cocycle  at  the element $a^p$. This
means that $\cry{C_{p^2}; X_i; T_i}$ is  torsion free.
\enddemo
\rightline{\text{\qed}}

\bigskip

Let $Y_i=\gp{\,\,\, \Phi(a),\,\, (a-1)^i\,\,\,}$  be a
$KC_{p^2}$-submodule in $KC_{p^2}$, where
$i=0, 1, \ldots, p-1$.
 The $K$-representation $\Gamma_i'$ corresponding to $Y_i$ has the
following form
$$
\Gamma_i'(a)=
\left(\smallmatrix
1 &   \gp{1}_0          & 0   \\
  &   \delta_1(a)       &  \gp{\alpha_i}_1      \\
  &                     &  \delta_2(a)          \\
\endsmallmatrix\right), $$
where $\alpha_i=(\xi_1-1)^i$ and $i=0, 1,\ldots, p-1$.

\proclaim
{Lemma 9}
For an arbitrary  cocycle $T:C_{p^2}\to
\widehat{Y_i}=FY_i^+/Y_i^+$  the group $\cry{C_{p^2}; Y_i; T}$
contains an element of order $p$.
\endproclaim
\demo{Proof} It is easy to see that an arbitrary cocycle of
$C_{p^2}$ with a value in $\widehat{Y_i}$ will be cohomologous to
a cocycle $T$, such that $T(a) =\lambda p^{-2}u+Y_i$, where
$\lambda\in K$, $u=\Phi(a)\Phi(b)$.   Thus, $T(a^p)=pT(a)=\lambda
p^{-1}u+Y_i$,  so to prove   the lemma it is sufficient to show
that $p^{-1}u\in(a^p-1)FY_i+Y_i$.  It is easy to see that
$$
(a^{p-1}+a^{p-2}+\cdots+a+1)-(a-1)^{p-1}=p\omega_1(a), \tag2
$$
where $\omega_1(a)\in KC_{p^2}$.

Let $v_1=(a-1)^i$.  Then from (2) it follows that
$$
(\Phi(a^p)-p)(a-1)^{p-i-1}v_1=
u-p\omega_1(a)\Phi(a^p)-p(a-1)^{p-i-1}v=u+py,
$$
where $y\in Y_i$.   Since $\Phi(a^p)-p=(a^p-1)z$, where
$z\in KC_{p^2}$, we have that $p\m1u+Y_i=
(a^p-1)p\m1z+Y_i$ which completes  the proof of the lemma.
\enddemo
\rightline{\text{\qed}}

\bigskip

Let $p\not=2$.  In the free $KC_{p^2}$-module $ (KC_{p^2})^{(2)}=
\{ \,(x_1,x_2) \, \mid \, x_1,x_2\in KC_{p^2} \, \} $ let us
consider $KC_{p^2}$-submodule
$$
U_j=\gp{ \,\,\,\,\, (\,
    (a-1)^{j+1}+\Phi(a),
\,
    (a-1)^{j}
\,
    );
\,\,\,\,\, \Phi(a^p)(\,a-1,\,1\,) \,\,\,\,\, },
$$
where $1\leq j\leq p-2$. The $K$-representation of $C_{p^2}$,
corresponding to the module $U_j$,  has the following form
$$
\Gamma_j'': a\to
\left(\smallmatrix
1 &        0            &    0           & \gp{1}_0            \\
  &        1            &    \gp{1}_0    & 0                   \\
  &                     &    \delta_1(a) & \gp{\alpha_j}_1     \\
  &                     &                & \delta_2(a)         \\
\endsmallmatrix\right),\tag3
$$
where $\alpha_j=(\xi_1-1)^j$ and $j=1,2,\ldots, p-2$.
Let us define a cocycle
$$
f_j:C_{p^2}\to \widehat{U_j}=FU_j^{+}/U_j^{+}
$$
by $f_j(a)=p^{-2}\Phi(a)\Phi(a^p)(1,0)+U_j$.

\proclaim
{Lemma 10}
The group $\cry{C_{p^2}; U_j; f_j}$ is  torsion free
( $j=1, \ldots, p-2$ ).
\endproclaim
\demo{Proof} Let $u_1=\Phi(a)\Phi(a^p)(1,0)$ and
$u_2=\Phi(a)\Phi(a^p)(0,1)$. It is easy to see that  the sequence
of $KC_{p^2}$-modules
$$
0\to Ku_2 \to U_j \to X_j\to 0\tag4
$$
is  exact.  The cocycle $f_j$ induces  a cocycle
$T_j:C_{p^2}\to \widehat{X_j}$ (see (1) which is not
equal to the zero cocycle  in the group $H=\gp{a^p}$ according
to  Lemma 8).   Therefore $f_j$ is also non-cohomologous to
the zero cocycle  in $H$.  This means that $\cry{C_{p^2};
U_j;  f_j}$ has no elements of order $p$.
\enddemo
\rightline{\text{\qed}}

\bigskip
Let us consider one more  $KC_{p^2}$-module:
$U_0=KC_{p^2}\Phi(a)$ generated by $\Phi(a)$ in
$KC_{p^2}$.  The $K$-representation of the group $C_{p^2}$
corresponding  to this module has the following form:
$$
a\to
\left(\smallmatrix
1 &   \gp{1}_0         \\
0 &  \delta_2(a)          \\
\endsmallmatrix\right). $$

\proclaim
{Lemma 11}
For any cocycle $T:C_{p^2}\to U_0$  the group
$\cry{C_{p^2};  U_0;  T}$ contains an element of order $p$.
\endproclaim
\demo{Proof}
It is easy to see  that any $1$-cocycle of the group
$C_{p^2}$ with values   in the group
$\widehat{U_0}=FU_0^{+}/U_0^{+}$ is cohomologous to the
cocycle $T$,  such that
$$
T(a)=\lambda p^{-2}\Phi(a)\Phi(a^p)+U_0,
$$
where $\lambda \in K$.
Replacing
   $a$ with  $a^p$ in (2)
we have
$$
{\textstyle\frac{1}{p}}\Phi(a)\Phi(a^p)=
{\textstyle\frac{1}{p}}(a^p-1)^{p-1}\Phi(a)+
\omega_1(a^p)\Phi(a).
$$
Then $ T(a^p)=(a-1)z+U_0, $ where $z\in FU_0$, which is  proves
the lemma.
\enddemo
\rightline{\text{\qed}}

\subhead
 Proof of  Theorem 2
\endsubhead
From the description of the $K$-representations  of the group
$C_{p^2}$ \cite{2},  it follows that all indecomposable
$KC_{p^2}$-modules of those faithful $K$-representations of the
group  $C_{p^2}$, whose  characters contain the trivial character
of $C_{p^2}$, are   the following:
$$
\split
  X_i \quad (i=0,1,\ldots, p-2); \quad
& Y_j \quad (j=0,1,\ldots, p-1); \quad  U_0; \\
& U_k \quad (k=1,\ldots,p-2). \\
\endsplit
$$

According to  Lemmas 9 and 11  we are interested  only in the
modules $X_i$ and $U_j$. Let us consider the module $X_i$ ($0\leq
i\leq p-2$). It is easy to see that  Lemma 2 can be applied to the
factor module $X_i/Kv$, where $v=\Phi(a)\Phi(a^p)$.  Therefore,
any cocycle of the group $C_{p^2}$ with the values in
$\widehat{X_i}$ will be cohomologous to such a  cocycle $T$ that
$$
T(a)=\lambda p^{-2}v+X_i,\tag5
$$
where $\lambda\in K$.  We will show that if in this equation
$\lambda\equiv 0{\pmod p}$  then the cocycle $T$ is cohomologous
to the trivial cocycle.  From (2) it follows that
$$
{\textstyle\frac{1}{p}}\Phi(a)\Phi(a^p)+
{\textstyle\frac{1}{p}}\Phi(a^p)(a-1)^{i+1}=
{\textstyle\frac{1}{p}}(a^p-1)^{p-1}\theta_i+\omega_1(a^p)\theta_i,\tag6
$$
where $\theta_i=\Phi(a)+(a-1)^{i+1}\in X_i$.
We will use the equation
$$
\Phi(a^p)(a-1)^{p}={p}(a-1)\Phi(a^p)\omega_2,
$$
where $\omega_2$ is a unit in $KC_{p^2}$.
From (6) it follows that
$$
{\textstyle\frac{1}{p}}(a-1)^{i+1}\Phi(a^p)=
{\textstyle\frac{1}{p^2}}\Phi(a^p)(a-1)^{p+i}\omega_2\m1\in (a-1)FX_i,
$$
for all $i=0,1,\ldots,p-2$.  Then from (6) one finds
$$
{\textstyle\frac{1}{p}}\Phi(a)\Phi(a^p)\in (a-1)FX_i+X_i
$$
for all $i=0,1,\ldots,p-2$.
This means that if in (5) $\lambda\equiv 0{\pmod p}$, then
the cocycle $T$ is cohomologous to the zero cocycle.

>From the above it follows that $H^{1}(C_{p^2},\widehat{X_i})$
is a cyclic group of order $p$ and all elements of this group can
be represented by the  cocycles $T$ of (5) with  $\lambda=0,1,
\ldots, p-1$.

We will show that each nonzero cocycle  $T$ defines up to
isomorphism one group $\cry{C_{p^2}; X_i; T}$.

Let $\varepsilon$ be an automorphism of the group
$C_{p^2}$ and $X_i^{\varepsilon}$ be the
$KC_{p^2}$-module $X_i$ twisted by this automorphism, i.e.
$$
X_i^\varepsilon=X_i, \,\,\,\, \,\,\,\, a\cdot x= \varepsilon(a)x,
\,\,\,\, \,\,\,\, x\in X_i.
$$
It is not difficult to show the existence of an  automorphism
$\tau$ of the $K$-module $X_i$ such  that $\varepsilon(a)\tau=\tau
a$ in $X_i$ and $\tau(v)=v$, where $v=\Phi(a)\Phi(a^p)$.

Let $\varepsilon\m1(a)=a^s$,   with  $(s,p)=1$.  Since
$aT_i(a)=T_i(a)$ and $\tau'(\ov1{v})=\ov1{v}$, where $\ov1{v}=v+X_i$,
we have that
$$
T_i^\varepsilon(a)=\tau'T_i(\varepsilon\m1(a))=sT_i(a)=
sp^{-2}v+X_i.
$$
From Lemma 1 it follows that $\cry{C_{p^2};  X_i;  T_i}$ is
isomorphic to $\cry{C_{p^2};  X_i;  T}$, where
$T(a)=sp^{-2}v+X_i$.  We have  shown that all groups
$\cry{C_{p^2};  X_i;  T}$, where $T\not\equiv 0$,  are
isomorphic to $\cry{C_{p^2};  X_i;  T_i}$.

Now let us consider the group $\cry{C_{p^2};  U_j;  T}$.
First of all we  remark that the group
$H^1(C_{p^2},\widehat{Y_j})$ is a cyclic
group of order $p$, where  $j=1,\ldots, p-1$.
The proof of this fact  is similar  to
the proof for  the group $H^1(C_{p^2},\widehat{X_i})$.  Since
$Y_0=KC_{p^2}$, we have that
$H^1(C_{p^2},\widehat{Y_0})=0$ (see Lemma 3).

Let $u_1,u_2,\ldots,u_{p^2+1}$ be a $K$-basis in $U_j$,
such that
$$
u_1=\Phi(a)\Phi(a^p)(1,0)
\,\,\,\,\,\,
\text{  and  }
\,\,\,\,\,\,
u_2=\Phi(a)\Phi(a^p)(0,1).
$$
We will use the exact sequence  (4) and the exact sequence
$$
0@> >> Ku_1 @> >>  U_j @> >>  Y_j@> >>  0.  \tag7
$$
This will give us an opportunity to show that any cocycle
$T:C_{p^2}\to U_j$ is cohomologous to a  cocycle
$T_{\alpha,\beta}$, such that
$$
T_{\alpha,\beta}(a)=p^{-2}(\alpha u_1+\beta u_2)+U_j,
$$
where  $0\leq \alpha,\beta\leq p-1$.

If $\alpha=0$ then according to  Lemma 9 and (7) it follows that
the cocycle $T_{0,\beta}$ is cohomologous to the zero cocycle  at
the element $a^p$ of the group $C_{p^2}$ and, therefore,  the
group $\cry{C_{p^2}; U_j; T_{0,\beta}}$ will  contain an element
of order $p$.

Now let  $\alpha\not=0$.   Then $\alpha$ is a unit in $K$ and
$\tau(x)=\alpha x$  ($x\in U_j$)   is an automorphism of
$KC_{p^2}$-module $U_j$.   From this it follows that the cocycle
$T_{\alpha,\beta}$ ($\alpha\not=0$) can be replaced  by
$T_{1,\alpha\m1\beta}$.  So  it is enough to consider the cocycles
$T_{1,\beta}$, where $\beta=0,1,\ldots,p-1$. We will show that
$\cry{C_{p^2}; U_j; T_{1,\beta}}$ is isomorphic to $\cry{C_{p^2};
U_j; f_j}$ (note that $f_j=T_{1,0})$.

Let us replace the basis element $u_1$ by   $u_1'=u_1+\beta
u_2$  in $U_j$.   Then
$$
T_{1,\beta}(a)=p^{-2}u_1'+U_j.
$$
Let $Y_j'=U_j/Ku_1'$. Then  the
$K$-representation $\Gamma_j'''$ corresponding to  the
$KC_{p^2}$-module $Y_j'$ is
$$
\Gamma_j''': a\to
\left(\smallmatrix
1 &   \gp{1}_0          & \gp{-\beta}_0        \\
  & \delta_1(a)         &  \gp{\alpha_j}_1     \\
  &                     &  \delta_2(a)         \\
\endsmallmatrix\right).
$$
This representation is equivalent to $\Gamma_j'$. According
to  this equivalence we will replace the basis elements
$u_2,\ldots,u_{p^2+1}$ with  $u_2',\ldots,u_{p^2+1}'$. Then in the
$K$-basis $u_1', u_2',\ldots,u_{p^2+1}'$ the operator $a$
have the same matrix (3) as  in the basis $u_1, u_2, \ldots,
u_{p^2+1}$.  Let us define an  automorphism $\tau: U_j\to U_j$ of
the $K$-module $U_j$ by  $\tau(u_i')=u_i$, where $i=1,\ldots,
p^2+1$.  As it was shown, $\tau a=a\tau$. Moreover,
$$
\tau'T_{1,\beta}(a)=\tau'(p^{-2}u_1'+U_j)=p^{-2}u_1+U_j=f_j(a).
$$
The  isomorphism of the groups $\cry{C_{p^2};  U_j;  T_{1,\beta}}$
and $\cry{C_{p^2};  U_j;  f_j}$  follows from  Lemma 1. So  we
have  shown that among the  groups $\cry{C_{p^2};  M; T}$ only
those can be indecomposable and  torsion free, for whose modules
$M$ and the cocycles $T$ were characterized  in this theorem.

Then Lemmas 8 and 10 complete  the proof.

\rightline{\text{\qed}}

Let $G\cong C_p \times C_p$ with  generators $a,b$ and let $K$ be
either   the ring of rational integers $\Bbb Z$, or the ring of
$p$-adic integers $\Bbb Z_{p}$,  or the localization  $\Bbb
Z_{(p)}$  of $\Bbb Z$ at $p$, respectively.   In case $p=2$ we
will give full description of the indecomposable  torsion free
groups $\cry{C_2 \times C_2; M; T}$.

In this case we will use the classification of the indecomposable
$K$-representations of the group $C_2 \times C_2$, given by
L.~Nazarova  in \cite{14, 15}.

\proclaim
{Lemma 12 }
Let $M$ be a $K[C_p \times C_p]$-submodule in the free
$K[C_p \times C_p]$-module $(K[C_p\times C_p])^{(2)}$ with
the following system of generators:
$$
M=\gp{
\,\,\,
(\Phi(a),0);
\,\,\,
(p,0);
\,\,\,
(0,\Phi(b));
\,\,\,
(0,p);
\,\,\,
(b-1,1-a)
\,\,\,}.
$$
Then the following conditions are  satisfied:
\itemitem{1)} $M$ is an indecomposable
$K[C_p \times C_p]$-module and $dim_K(M)=2p^2$;
\itemitem{2)} there exists a cocycle $T:C_p\times C_p\to
\widehat{M}=FM^{+}/M^{+}$, such that:
$$
T(a)=(1,0)+M,
\,\,\,\,\,\,\,\,\,\,\,\,\,
T(b)=(0,1)+M;
$$
\itemitem{3)} the group $\cry{C_p \times C_p; M; T}$
is torsion free.
\endproclaim
\demo{Proof}
Let $\ov1{Z_p}=K/pK$.  Obviously, $\ov1{Z_p}$ is a
$K[C_p\times C_p]$-module with  trivial action  of   $C_p \times
C_p$.  Let us consider the projective resolution  of $K[C_p
\times C_p]$-module $\ov1{Z_p}$:
$$
\CD
\cdots @>{ }>>
(K[C_p\times C_p])^{(3)}
@>\tau_1 >>
(K[C_p\times C_p])
@>\tau_0 >>
\overline{Z_p}
\to 0.
\\
\endCD\tag8
$$
It is  easy to see that
$\ker(\tau_0)=\gp{
\,\,\,\,\,
a-1,
\,\,\,\,\,\,
b-1,
\,\,\,\,\,
p
\,\,\,\,\,
}$, and
$$
\split
\ker(\tau_1)=\gp{
\,\,\,
(\Phi(a),0,0);
\,\,\,\,
(0,\Phi(b),0); &\\
\,\,\,\,
(b-1,1-a,0);   &
\,\,\,\,
(p,0,1-a);
 \,\,\,\,
(0,p,1-b)
\,\,\,
 }. \\
\endsplit
$$
The $K[C_p \times C_p]$-modules $\ker(\tau_0)$ and $\ker(\tau_1)$
are indecomposables   of the $K$-representations of $C_p \times
C_p$.  Each  $x\in \ker(\tau_1)$  has the following form:
$$
x=(\,\,\,
   u_1\Phi(a)+u_3(b-1)+pu_4; \,\,\,
   u_2\Phi(b)+u_3(1-a)+pu_5; \,\,\,
   u_4(1-a)+u_5(1-b)\,\,\,),\tag9
$$
where $u_i\in K[C_p \times C_p]$ and $i=1,\ldots,5$.

We will assign to each    $x\in \ker(\tau_1)$ ( see (9)) the element
$$
(\,\,\, u_1\Phi(a)+u_3(b-1)+pu_4; \,\,\,
   u_2\Phi(b)+u_3(1-a)+pu_5\,\,\,)
$$
of $M$.  It is easy to check that this map generates  an
isomorphism of the $K[C_p \times C_p]$-modules $\ker(\tau_1)$ and
$M$. Thus, we have shown  that $M$ is an indecomposable  $K[C_p
\times C_p]$-module. Since $dim_K(T_0)=p^2$, we have
$$
dim_K(M)=dim_K(\ker(\tau_1))=dim_K(K[C_p \times C_p])^{(3)}
-dim_K(\ker(\tau_0))=2p^2.
$$

2) Let us define  a  function
$T:C_p \times C_p\to \widehat{M}$ by:
$$
T(a^i)=(\,\,\, 1+a+\cdots+a^{i-1}, \,\,\,0\,\,\, )+M;
$$$$
T(b^j)=(\,\,\,0,\,\,\,1+b+\cdots+b^{j-1}\,\,\,)+M;
$$$$
T(a^ib^j)=a^{i}T(b^j)+T(a^i)+M;
\,\,\,\,\,
T(1)=M,
$$
where $0< i,j\leq p-2$.  It is easy to see that
$\Phi(a)T(a)\subset M$,
$\Phi(b)T(b)\subset M$ and
$(a-1)T(b)-(b-1)T(b)\subset M$.
It follows that $T$ is a $1$-cocycle of
$C_p \times C_p$ with values in $\widehat{M}=FM^{+}/M^{+}$.

3) It is sufficient to show that the cocycle $T$ is not
cohomologous  to  the zero cocycle  at every nontrivial element
$g$ of $C_p\times C_p$.  Let $g=a^ib^j$, where  $0< i,j\leq p-1$.
Suppose that there  exists  $z\in FM$, such that
$$
T(g)=(g-1)z+M. \tag10
$$
Then from the definition of $T$ and from (10) it follows that
$$
(\,\,\,1+a+\cdots+a^{i-1},\,\,\,a^i(1+b+
\cdots+b^{j-1})\,\,\,)=(g-1)z+x,
$$
where $x\in M$.
Multiplying the last equality  by $\Phi(a)\Phi(b)$ taking into account  that
$$
\Phi(a)\Phi(b)M=p\Phi(a)\Phi(b)(K,K),
$$$$
\Phi(a)\Phi(b)(g-1)=0
$$
can conclude  that $(i,j)\in (pK,pK)$, which is impossible for
$0<i,j\leq p-1$.  This  contradicts  to the
assumption that the cocycle $T$ is cohomologous   to the zero
cocycle  at  $g$ (see (10)).

Similarly, we show that the cocycle $T$ is not  cohomologous
to the zero cocycle at the  rest of the nontrivial elements of the
group $G\cong C_p\times C_p$.  Thus, the group
$\cry{C_p\times C_p; M; f}$ is torsion-free.
\enddemo\rightline{\text{\qed}}

Let $p=2$, $G=\gp{a}\times\gp{b}\cong C_2\times C_2$ and let $K$
be either   the ring of rational integers $\Bbb Z$, or the ring of
$2$-adic integers $\Bbb Z_{2}$,  or the localization  $\Bbb
Z_{(2)}$  of $\Bbb Z$ at $2$, respectively. We will study those
groups  $\cry{G; M; T}$ which are torsion free.

The group  $G$ has the following irreducible
$K$-representations:
$$
\matrix
\format
\quad \l\quad
&\quad \l  &\r
&\quad \l  &\r
&&
\quad \l\quad &\quad \l\quad &\quad \l &\r
&\quad \l &\r \\
\chi_0: & a\to &1, &  b\to &1; &    & \chi_1:& a\to &-1,& b\to &1;  \\
\chi_{2}: & a\to &-1,&  b\to &-1;&    & \chi_3:& a\to &1, & b\to &-1.  \\
\endmatrix
$$
Let $H=\gp{h}$ be a subgroup of $G$ of order $2$. The
indecomposable $K$-representations of $H$, up to equivalence, are
one of the following:
$$
\gamma_0:  h \to 1;
\,\,\,\,\,\,
\gamma_1:  h \to -1;
\,\,\,\,\,\,
\gamma_2:  h \to \left(\smallmatrix
                   1 &  1\\
                   0 & -1\\
                 \endsmallmatrix\right).
\tag11
$$
Let $\Gamma$ be a $K$-representation of $G$, $\Gamma|_{H}$ be the
restriction of  $\Gamma$ to the subgroup $H$, $M$ be a $KG$-module
of the $K$-representation $\Gamma$ and $T:G\to \widehat{M}$ be an
arbitrary  cocycle of $G$ with values in the group
$\widehat{M}=FM^+/M^+$ ($F$ is a field containing $K$).  The
following Lemma gives the necessary conditions   for $\cry{G;M;T}$
to  be torsion free.

\proclaim {Lemma 13} If $\cry{G; M; T}$ is torsion free then for
any one of the three nontrivial elements non-trivial subgroups $H$
of order $2$, the restriction $\Gamma|_{H}$ of the representation
$\Gamma$ contains the trivial representation $\gamma_0$ in the
decomposition of $\Gamma|_{H}$ into the direct sum of
indecomposable $K$-representations  of $H=\gp{h}$.
\endproclaim

\demo{Indirect proof} Assume that there is  a nontrivial element
$h$ in $G$, such that the $K$ representation $\Gamma|_{h}$ of $H$
is the direct sum  of $K$-representation $\gamma_1$ and $\gamma_2$
( with the multiples ) but does  not include $\gamma_0$ in this
sum (see (11)).  Then from  Lemma 2 and 3 it follows that any
cocycle $T:G\to \widehat{M}$ in  the subgroup $H=\gp{h}$ will be
cohomologous to the zero cocycle, which implies the  existence of
elements of  order $2$ in  $\cry{G; M; T}$.
\enddemo
\rightline{\text{\qed}}

\bigskip

We make some remarks about the $K$-representations of
$G\cong C_2\times C_2$.  Let $K$ be a $KG$-module of the trivial
representation $\chi_0$ of $G$.   Let us consider the
projective resolution of  the $KG$-module  $K$:
$$
\CD
\cdots \to
(KG)^{(n)}
@>\nu_n >>
(KG)^{(n-1)}
@> >>
\cdots
\\
\cdots
@>\nu_3 >>
(KG)^{(2)}
@>\nu_2 >>
(KG)
@>\nu_1 >>
K
\to 0,
\\
\endCD\tag12
$$
where $\nu_n$ is  a homomorphism of the $KG$-modules
($n=1,2,\ldots $).

Then $\ker(\nu_n)$ is an indecomposable $KG$-module, and
$$
dim_K(\ker(\nu_n))=2n+1,
$$
where $n=1,2,\ldots$.

Let $\Gamma_n$ be the $K$-representation of $G$ corresponding to
some $K$-basis in $\ker(\nu_n)$, and let $\Gamma_n^{*}$ be the
contragradient $K$-representation of $\Gamma_n$, that is
$\Gamma_n^*(g)=\Gamma^{\scriptstyle{T}}(g\m1)$, for all $g\in G$.

\proclaim {Lemma 14} (see \cite{15, 21}) All indecomposable,
pairwise nonequivalent $K$-represen\-tations  of $G\cong C_2\times
C_2$ of odd degree  are either $\chi_i$ or the tensor product
$\Gamma_n\otimes_K\chi_i$ or $\Gamma_n^*\otimes_K\chi_i$, where
$i=0, 1, 2, 3$ and $n=1, 2, \ldots $.
\endproclaim
\rightline{\text{\qed}}

\bigskip

Let $p=2$ in (8) and let us consider the projective
resolution  for
$\ker(\tau_0)=\gp{\,\, a-1,\,\, b-1,\,\, 2\,\,\,}$:
$$
\CD
\cdots \to
(KG)^{(t_n)}
@>\tau_n >>
(KG)^{(t_{n-1})}
@> >>
\cdots
\\
\cdots
@>\tau_3 >>
(KG)^{(t_2)}
@>\tau_2 >>
(KG)^{(t_1)}
@>\tau_1 >>
\ker( \tau_0)
\to  0.
\\
\endCD\tag13
$$
It is easy to show that in (13)
$t_n=2n+1$ and
$$
dim_K(\ker(\tau_n))=4n+4,
$$
where $n\geq 0$. Moreover,  all $KG$-modules $\ker(\nu_n)$ are
indecomposable.   If we take the tensor product of the
exact sequence
(12) over the  ring $K$ by the $KG$-module $\ker(\tau_0)$ and
compare  the result with the sequence (13), then it is easy to get
the isomorphism:
$$
\ker(\tau_0) \otimes_K \ker(\nu_n)\cong
\ker(\tau_n)\oplus P_n,
$$
where $P_n$ is a projective $KG$-module.

\proclaim
{Lemma 15}
Let $W_n$ be a  $K$-representation of
$G=\gp{a}\times\gp{b}\cong C_2\times C_2$
corresponding to the module $\ker(\tau_n)$ ($n\geq 0$).
This representation has the following form:
$$
\matrix
&&&&\\
W_0:&
a\to &
\left( \smallmatrix
1  &  1    &  0   &  1    \\
   & -1    &  0   &  0    \\
   &       &  1   &  0    \\
   &       &      & -1    \\
\endsmallmatrix\right),
&
b\to &
\left( \smallmatrix
1  &  1    &  1   &  0    \\
   & -1    &  0   &  0    \\
   &       & -1   &  0    \\
   &       &      &  1    \\
\endsmallmatrix\right); \\ &&&&\\
&&&&\\
W_n:&
a\to&
\left(\smallmatrix D        & 0   & 0    &  0       &  0        \\
          & E_n & 0    &  0       &  V_n      \\
          &     & -E_n &  V_n     &  0        \\
          &     &      &  E_{n+1} &  0        \\
          &     &      &          & -E_{n+1}  \\
\endsmallmatrix\right),  &
b\to&
\left(\smallmatrix
D         &  0  &  0   &     S       & 0         \\
          & E_n &  0   &   {V_n}'    & 0         \\
          &     & -E_n &     0       & V_n'       \\
          &     &      &  -E_{n+1}   & 0         \\
          &     &      &             & E_{n+1}   \\
\endsmallmatrix\right),\\  &&&&\\
\endmatrix
$$
where $D=\left(\smallmatrix
1  &  1    \\
0  & -1    \\
\endsmallmatrix\right)
$,
$S=\left(
  \smallmatrix
0  &\ldots 0  &  1 & 1    \\
0  &\ldots 0  &  0 & 0    \\
  \endsmallmatrix
  \right ),
$ and
$V_n=(0 E_n)$, ${V_n}'=(E_n 0)$
are matrices with $n$ rows and $n+1$ columns
$(n\geq 1)$.
\endproclaim

\demo{Proof} The proof of the  Lemma reduces  to the determination
of a $K$-bases of $ker(\tau_n)$, which is not difficult  to
construct inductively with respect to  $n$.
\enddemo
\rightline{\text{\qed}}

\proclaim
{Lemma 16}
A faithful  indecomposable  $K$-representation of
$G=\gp{a}\times\gp{b}\cong C_2\times C_2$ which
satisfies the necessary condition for the existence  of the torsion-
free group  $\cry{G; M; f}$ is one of the following:
$$
\matrix
           &&&&&&&&&\\
\Delta_n:  &&& a \to \Delta_n(a),
           &&& b \to \Delta_n(b),
           &&&   (n\geq 1);\\
           &&&&&&&&&\\
\Delta_n^*:&&& a \to \Delta_n(a)^{\scriptstyle{T}},
           &&& b \to \Delta_n(b)^{\scriptstyle{T}},
           &&&   (n\geq 1);\\
           &&&&&&&&&\\
W_n:       &&& a \to W_n(a),
           &&& b \to W_n(b),
           &&&   (n\geq 0);\\
           &&&&&&&&&\\
W_n^*:     &&& a \to W_n(a)^{\scriptstyle{T}},
           &&& b \to W_n(b)^{\scriptstyle{T}},
           &&&   (n\geq 0),\\
           &&&&&&&&&\\
\endmatrix
$$
where $A^{\scriptstyle{T}}$ is the transpose of $A$;
$$
\Delta_n(a)=
\left(\smallmatrix
E_n  & 0 &  0   &  E_n & 0    \\
     & 1 &  0   &  0   & 0    \\
     &   & -E_n &  0   & E_n  \\
     &   &      & -E_n & 0    \\
     &   &      &      & E_n  \\
\endsmallmatrix\right),
\quad
\quad
\quad
\quad
\Delta_n(b)=
\left(\smallmatrix
1    &   0 &  0   &  0   &  0    \\
     & E_n &  0   &  0   &  E_n  \\
     &     & -E_n &  E_n &  0    \\
     &     &      &  E_n &  0    \\
     &     &      &      & -E_n  \\
\endsmallmatrix\right);
$$
$\Delta_n^*$  and   $W_n^*$ are $K$-representations
of $G$,   contragradient to
the $K$-representations
$\Delta_n$ and $W_n$, respectively.
\endproclaim

\demo{Proof}
All $K$-representations listed above satisfy the necessary
condition for the existence of a torsion free group
$\cry{G; M; f}$.  The analysis  of all
representations of odd degree (see Lemma 14) shows that
among the representations $\Gamma_n\otimes \chi_i$ the
necessary condition is satisfied only by $\Delta^*$ which
is equivalent to $\Gamma_{2n}$ $(n=1,2,\ldots$).   Besides,
the representations
$W_n$ and $W_n^*$  the group $G$  has a parameterized series
of representations, whose degrees are divisible by $4$.
In this series the following pairs of matrices
correspond to the pairs of generator elements of $G$:
$$
\left(\smallmatrix
E_n  &  0   & 0    &  E_n \\
     & -E_n & E_n  &  0   \\
     &      & E_n  & 0    \\
     &      &      & -E_n   \\
\endsmallmatrix\right),
\quad\quad
\left(\smallmatrix
E_n &  0  &  \frak{F}&  0   \\
    &-E_n &    0     &  E_n \\
    &     & -E_n     &  0   \\
    &     &          &  E_n \\
\endsmallmatrix\right), $$
where the $K$-matrix $\frak{F}$ has  indecomposable  modulo $2K$
the Frobenius (i.e. rational) canonical normal form. Clearly, the
representation of this series does not satisfy the necessary
condition.  Let us consider the following pair of matrices:
$$
\left(\smallmatrix
1 & 0  &  0   &  0    & 0   & 0   \\
  & E_n&  0   &  0    & 0   & E_n \\
  &    & -E_n & E_n   & 0   & 0   \\
  &    &      & E_n   & 0   & 0   \\
  &    &      &       & 1   & 0   \\
  &    &      &       &     & -E_n \\
\endsmallmatrix\right),
\quad\quad
\left(\smallmatrix
E_{n+1} & 0   & E_{n+1} &  0  \\
        &-E_n &  0      &  E_n\\
        &     & -E_{n+1}&  0  \\
        &     &         &  E_n\\
\endsmallmatrix\right).
$$
These matrices define the indecomposable $K$-representation of $G$
of degree  congruent to $2$  modulo $4$. Obviously,  this
representation does not satisfy the necessary condition. We can
obtain the remaining representations of  degree $4n-2$ either by
the presented process of tensor multiplication by irreducible
$K$-representation  or by taking the contragradient
representation.

As a result we  get the representations which do not
satisfy the necessary condition for the existence of the
torsion free group $\cry{G; M; f}$.   Thus,  we
have considered  all  indecomposable $K$-representations  of
$G$.  The lemma is proved.
\enddemo
\rightline{\text{\qed}}

\subhead
Proof of  Theorem 3
\endsubhead
Note that, the module  $M$ of a $K$-representation $\Gamma$ of $G$
is  a module with  $m$-dimensional columns, whose entries
belong to the ring $K$ ( $m$ is the degree of the  representation
$\Gamma$).  Then $FM$ is a linear space of $m$-dimensional vector
columns over the field $F$, $\widehat{M}=FM^+/M^+$ is a group of
$m$-dimensional columns, whose  coordinates  belong to the group
$\widehat{F}=F^+/K^+$. Let $f:G\to \widehat{M}$ be a cocycle.  The
value $f(g)$ of $f$ at $g\in G$ is an $m$-dimensional column over
$\widehat{F}$. We note that if $g,h\in G$ then the multiplication
$g\cdot f(h)$ is the multiplication of the matrices $\Gamma(g)$
and $f(h)$.

If we consider the  coordinates of the
vector $f(g)$ as  elements of the field $F$, then
the elements of the ring $K$ will be replaced by $0$.

Let $\Gamma$ be any one  of the representations of $G$ listed in
Lemma 16, and let $M$ be the module of this representation, and
$H=\gp{h}$ a non-trivial subgroup of $G$. There exists only one
basis vector $v$ in $M$ such that the $KH$-module $M$ is a direct
sum $M=Kv\oplus M'$ of the $KH$-module $Kv$ and the $KH$-module
$M'$, generated by the rest of the basis vectors of $M$.   In
addition $hv=v$ and a $K$-representation $\Gamma'$ of $H$
corresponding to the module $M'$ is a sum of representations
$\gamma_1$ and $\gamma_2$ (see (11) ).  This allows us the
possibility to replace the  cocycle $f$ by  the cohomologous
cocycle $f_1$ in such a way that the projection
${f_1}_{|_{\widehat{M'}}}$ will be equal to zero in the element
$h$ (see Lemmas 2--3).

The coordinate   $x_v$ of the vector $f(h)$,   corresponding to
the  basis vector  $v$,  will be called   the {\it special
component} of the vector $f(h)$.   From $(1+h)f(h)=0$ it follows
that $2x_v=0$ (in the group $\widehat{F}$).  Besides, for any
vector $z\in \widehat{M}$  the special component of the vector
$(h-1)z+f(h)$ is always equal to $x_v$.  If $x_v=\frac{1}{2}$, then the
cocycle $f$ is  not cohomologous to  the zero cocycle  at  $h$.

These remarks justify the following plan for the
construction of cocycles of the representations $\Gamma$
from Lemma 16.  The form of the representation $\Gamma$
 defines  the special components of the vectors $f(a)$ and
$f(b)$ ( $a$ and $b$ are generators of $G$).   We choose the
vector $f(a)$  such that we deduce that the special component is
$\frac{1}{2}$ and all the rest are zero.
The possible forms of  the components of the vector
$f(b)$ follow from  the conditions:
$$
(1+b)f(b)=0,\tag14
$$$$
(a-1)f(b)=(b-1)f(a).  \tag15
$$
We  will carry out the following  operations  on
the vector $f(b)$:   replace  $f(b)$ by  the vector
$$
(b-1)z+f(b),\tag16
$$
where $z\in\widehat{M}$ and $(a-1)z=0$.

We discard all those  forms of  $f(b)$, with a zero special component.
For the vector $f(b)$ whose special  component equals to
$\frac{1}{2}$,  we find:

$$
f(ab)=af(b)+f(a)\tag17
$$
and we  examine the solvability of the following equation
$$
(ab-1)z+f(ab)=0\tag18
$$
with the variable $z\in \widehat{M}$.

The group $\cry{G; M; f}$ is torsion free if and only if
 the equation (17) is unsolvable.

We  consider the following seven  cases:

Case 1.  Let  $\Gamma=\Delta_n$. The special components are the
$(n+1)$th entry in  $f(a)$ and the first one in  $f(b)$. Set the
$(n+1)$th coordinate  of  $f(a)$  to  $\frac{1}{2}$ and let all
the rest be  $0$. Let
$$
f^{\scriptstyle{T}}(b)=(y,Y_1,Y_2,Y_3,Y_4),\tag19
$$
where $y\in\widehat{F}$, $Y_i\in\widehat{F}^{(n)}$ and  $i=1,2,3,4$.

The  operation (16) can replace  $Y_2$ by  the zero vector. From
(15)  follows that $Y_3=Y_4=0$,  and from (14), it follows that
$2y=0$ and $2Y_1=0$.   Let $y=\frac{1}{2}$,
$Y_1=(v_1,v_2,\ldots,v_n)$.  Using  (17),  it is easy to transform
(18) to a linear  system of equations (over $\widehat{F}$) with
the $(n+1)\times n$-matrix
$$
\left(\smallmatrix
1  & 0 & \ldots & 0 & 0\\
-1 & 1 & \ldots & 0 & 0\\
 \ldots & \ldots & \ldots & \ldots & \ldots \\
 \ldots & \ldots & \ldots & \ldots & \ldots \\
0 & 0 & \ldots & -1 & 1\\
0 & 0 & \ldots & 0 & -1\\
\endsmallmatrix\right)
$$
and  free  coefficients   $\frac{1}{2}$, $v_1,\ldots,
v_{n-1}, v_n+\frac{1}{2}$.   This system is  solvable  if
and only if   $v_1+\cdots+ v_{n-1}+ v_n=0$.

Case 2.  Let $\Gamma=\Delta_n^*$.  The matrices  of the
$K$-representation were transposed   to the matrices of
$\Delta_n$.  The special components of  $f(a)$ and $f(b)$ are
the same as in Case 1.   Let
us assume that $f(a)$ and $f(b)$ are   chosen at  first in
the same way  as in the case of $\Gamma=\Delta_n$ (see
(19)).   Condition (14) and operation (16) transform  the
vector $f(b)$ to the following form:
$$
f^{\scriptstyle{T}}(b)=(y,0,-2Y_3,Y_3,0).
$$
Let $Y_3=(v_1,\ldots, v_{n-1}, v_n)$.  It follows from (15) that
if $n\geq 2$ then
$$
\split
& y-2v_1=0;\,\,\,\,\,  2v_2=\cdots=2v_n=0;\\
& 2v_1=0; \ldots; 2v_{n-1}=0; \,\,\,\,\, 2v_n=
{\textstyle\frac{1}{2}}, \\
\endsplit
$$
and, if  $n=1$, then  $y-2v_1=0$, $2v_1=\frac{1}{2}$.

If $n>1$ and $y=\frac{1}{2}$,  then (19) leads to a contradiction.
If $n=1$ and $y=\frac{1}{2}$, then  $v_1=\frac{1}{4}$.

Thus, if $n>1$ and $f$ is a cocycle then the special
component of the vector $f(b)$ is equal to zero.  Then  the
cocycle $f$ is cohomologous  to  the zero cocycle  in the
element $b\in G$.  This  means that the  torsion free group
$\cry{G; M; f}$ does not  exist for the representation
$\Gamma=\Delta_n^*$, if  $n>1$.

Let $n=1$.  Then
$$
f(a)=(0,{\textstyle\frac{1}{2}},0,0,0);
\,\,\,\,\,
f(b)=({\textstyle\frac{1}{2}},0,{\textstyle\frac{1}{2}},
{\textstyle\frac{1}{4}},0);
\,\,\,\,\,
f(ab)=({\textstyle\frac{1}{2}},0,{\textstyle\frac{1}{2}},
{\textstyle\frac{1}{4}},{\textstyle\frac{1}{2}}).
$$
It is easy to check that (18) is  unsolvable.

Case 3.  Let $\Gamma=W_n^*$ ($n>0$).  The special components are
$(2n+3)$th of  $f(a)$ and the last of  $f(b)$.  Let the special
component of  $f(a)$ be equal to  $\frac{1}{2}$,  and let all the
rest be  zero.

Let $f^{\scriptstyle{T}}(b)=(Y_0, Y_1, Y_2, Y_3, Y_4)$, where $Y_0\in
\widehat{F}^{(2)}$, $Y_1, Y_2\in \widehat{F}^{(n)}$,
$Y_3, Y_4\in \widehat{F}^{(n+1)}$.   Operation (16)
allows us to replace $Y_3$  by the zero vector.  It follows
from (14) that $Y_1=0$.   Condition (15) shows that
$Y_2=0$, $Y_0=(0,y)$ ($y\in \widehat{F}$, $2y=0$) and $2Y_4=0$.
Consequently,
$$
f^{\scriptstyle{T}}(b)=(0, y, 0, \ldots, 0,v_1, \ldots, v_n,
{\textstyle\frac{1}{2}}).
$$
The special component of $f(ab)$ is the second coordinate which,
according to (17),  equals to $y$.  Therefore $y=\frac{1}{2}$ and for
any $v_1, \ldots, v_n$ ($2v_1=2v_2= \cdots=2 v_n=0$) the group
$\cry{G; M; f}$ is  torsion free.

Case 4.  Let $\Gamma=W_0^*$.  In this case it is easy to see  that the
cocycle $f$ with
$$
f(a)=(0, 0, {\textstyle\frac{1}{2}}, 0);
\,\,\,\,\,
f(b)=(0, {\textstyle\frac{1}{2}}, 0,{\textstyle\frac{1}{2}})
$$
determines  a  torsion free group $\cry{G; M; f}$.

Case 5.  Let $\Gamma=W_n$ ($n>1$).  We take the vectors $f(a)$ and
$f(b)$   in the same fashion as in the Case 3.
Condition (14) shows that all
components of the vector $Y_4$, except the last, are zero.
Then  condition  (15) leads  to   contradiction.

We obtain  the contradiction by setting the special component
in $f(a)$  equal to $\frac{1}{2}$.  Consequently, for
$\Gamma=W_n$ ($n>1$) any cocycle $f$ is cohomologous to
the zero cocycle at the generator $a$ of $G$, this means that
the torsion free group $\cry{G; M; T}$ does not exist in this
case.

Case 6.  Let $\Gamma=W_1$.  In a cocycle $f$ with
$$
f(a)=(0, 0, 0, 0, {\textstyle\frac{1}{2}}, 0, 0, 0);
\,\,\,\,\,
f(b)=(0, {\textstyle\frac{1}{2}}, 0,
{\textstyle\frac{1}{4}},0, {\textstyle\frac{1}{2}},
0,{\textstyle\frac{1}{2}})
$$
the special components of the vector $f(a)$ ( the fifth one
) and $f(b)$ ( the last one), and $f(ab)$ ( the second one)
are all equal to $\frac{1}{2}$.   The cocycle $f$ determines the
torsion free group $\cry{G; M; f}$ (see also the Lemma 12).

Case 7.  Let $\Gamma=W_0$.  The special components are the
$3^{rd}$ for $f(a)$ and the $4$th  for $f(b)$.  Let us assume that
they are equal to $\frac{1}{2}$.  Then there exists only one
cocycle
$$
f(a)=(0, 0, {\textstyle\frac{1}{2}}, 0);
\,\,\,\,\,
f(b)=(0, 0, 0,{\textstyle\frac{1}{2}}).
$$
Hence, $f(ab)=(\frac{1}{2}, 0, \frac{1}{2}, \frac{1}{2})$
and the special component  (the second one) for $f(ab)$
is equal to zero.   The cocycle $f$ is cohomologous to zero
on the element $ab$.  The group $\cry{G; M; f}$ contains
elements of  order $2$.

 It follows from Lemma 16 that all the $K$-representations
$\Gamma$ of $G$ have been enumerated, for  the torsion-
free group $\cry{G; M; f}$  to exist. Consequently, the Theorem is
proved.

\rightline{\text{\qed}}

\enddocument